\documentstyle{amsart}
\title{Dolbeault Cohomology of compact Nilmanifolds}
\author{Sergio Console \and Anna Fino}
\thanks{Research partially supported by MURST and CNR of Italy}
\address{
Dipartimento di Matematica\\ Universit\`a di Torino\\ Via Carlo 
Alberto 10\\ I-10123 Torino}
\email{console@@dm.unito.it, fino@@dm.unito.it}

\newcommand{\om}{\omega}

\newcommand{\beq}{\begin{equation}}
\newcommand{\eeq}{\end{equation}}
\newcommand{\bqn}{\begin{eqnarray}}
\newcommand{\eqn}{\end{eqnarray}}
\newcommand{\bqne}{\begin{eqnarray*}}
\newcommand{\eqne}{\end{eqnarray*}}
\newcommand{\cg} {{\frak g}}

\newcommand{\ct}{{\frak t}}

\newcommand{\ca}{{\frak a}}

\newcommand{\cb}{{\frak b}}

\newcommand{\ch}{{\frak h}}

\newcommand{\cf}{{\frak f}}

\newcommand{\R}{{\Bbb R}}
\newcommand{\Q}{{\Bbb Q}}

\newcommand{\C}{{\Bbb C}}
\newcommand{\de}{\partial}
\newcommand{\debar}{\overline \partial}
\newcommand{\w}{\wedge}

\newcommand{\HD}[2]{H_{\overline \partial}^{{#1},{#2}}}
\newcommand{\hd}{H_{\overline \partial}} 
\newcommand{\MM}{{\cal M}}

\newcommand{\lb}[1]{\label{#1}}
\newcommand{\rf}[1]{(\ref{#1})}\newcommand{\rl}[1]{\ref{#1}}

\newcommand{\ba}{\begin{array}}
\newcommand{\be}{\begin{equation}}
\newcommand{\ea}{\end{array}}
\newcommand{\ee}[1]{\label{#1}\end{equation}}

   at 14 truept  at 10 truept  at 12.00 truept

\begin{document}
\maketitle

\begin{abstract}
Let $M= G/\Gamma$ be a compact nilmanifold  endowed with an 
invariant complex
structure. We prove that, on an open set of any connected 
component of the moduli space ${\cal C} ({\frak g})$ of invariant 
complex structures on $M$, the Dolbeault cohomology of $M$ is 
isomorphic to the one of the
differential bigraded algebra associated 
to the complexification $\cg^\C$ of the Lie algebra of $G$.
To obtain this result, we first prove the above isomorphism for compact nilmanifolds endowed 
with a rational  
invariant complex structure.  This is done using a descending series associated to the complex 
structure and the Borel spectral
sequences for the corresponding set of holomorphic fibrations.
Then we apply the 
theory of Kodaira-Spencer for 
deformations of complex structures. 
\end{abstract}

\section{Introduction}\lb{Intr}

Let $M$ be a compact nilmanifold of real dimension $2n$. It follows 
from  a
result of Mal'\v{c}ev \cite{Ma} that $M=G/\Gamma$  where $G$ is a  
simply connected
$(s+1)$-step nilpotent Lie group admitting  a basis of left invariant 
1-forms
for which the coefficients in the structure equations are rational 
numbers,
and $\Gamma$ is a lattice in $G$ of maximal rank (i.e., a discrete 
uniform
subgroup, cf. \cite{Ra}).  We will let $\Gamma$ act on $G$ on the left. 
It is well known that such a lattice $\Gamma$ exists in $G$ if and 
only 
if the Lie algebra $\frak g$ of $G$  has a
rational structure, i.e. if there exists a rational Lie subalgebra 
${\cg}_{\Q}$ such that
${\cg} \cong  {\cg}_{\Q} \otimes \R$.

The de Rham
cohomology of a compact nilmanifold can be computed by means of the
cohomology of the Lie algebra of the corresponding nilpotent Lie group
(Nomizu's Theorem \cite{No}). 

We assume that $M$ has an invariant complex structure $J$, 
that is to say that
$J$ comes from a (left invariant) complex structure $J$ on $\cg$.
   
\noindent Our aim is to relate the Dolbeault cohomology of $M$ with 
the cohomology ring ${\HD {*} {*}} (\cg^\C)$  of the
differential bigraded algebra $\Lambda^{*,*}(\cg^\C)^*$, associated 
to 
$\cg^\C$ with respect to the operator $\debar$ in the canonical
decomposition $d=\de+\debar$ on $\Lambda^{*,*}(\cg^\C)^*$. 

The study of the Dolbeault cohomology of nilmanifolds with an 
invariant
complex structure is motivated by the fact that the latter provided 
the first
known examples of compact symplectic manifolds which do not admit any
K\"ahler structure \cite{Ab,CFG,Th}.

Since there exists a natural map
$$
i: {\HD {*} {*}} (\cg^\C) \rightarrow {\HD {*} {*}} (M)
$$
which is always injective (cf. Lemma 7), the problem we will study is to see for which complex 
structure $J$ on $M$ 
the above map gives an isomorphism
\begin{equation}
{\HD {p} {q}} (M) \cong {\HD {p} {q}} (\cg^\C). \label{iso}
\end{equation}
Note that $\HD {*} {*} (\cg^\C)$ can be identified with the 
cohomology of the
Dolbeault complex of the  forms on $G$ which are invariant by the 
left action
of $G$ (we shall call them  briefly $G$-invariant forms) and $\HD 
{*}{*} (M)$
with the cohomology of the Dolbeault complex of $\Gamma$-invariant 
forms on
$G$. We shall use these identifications throughout this note.

Our main result is the following

\medskip

\noindent{\bf Theorem A} 
{\sl The isomorphism \rf{iso} holds on an open set of any connected 
component of the moduli space ${\cal C} ({\frak g})$ of invariant 
complex structures on $M$.}

\medskip

To obtain Theorem A we first consider the case of complex structures 
$J$ which are rational, i.e.  they are
compatible with the rational structure of $G$ ($J 
({\cg}_{\Q}) \subseteq {\cg}_{\Q}$).

\medskip

\noindent{\bf Theorem B} 
{\sl For any rational complex structure $J$, the isomorphism \rf{iso}
holds. }

\medskip

It is an open problem whether the isomorphism (1) holds for any 
compact nilmanifold endowed with an arbitrary invariant complex
structure. We do not know examples for which (1) does not hold. 

\smallskip

Theorem A will follow from Theorem B using the theory of 
deformations of complex structures \cite{KS,Su}. 
Indeed by \cite{Sa} the set ${\cal C} ({\cg})$ of complex structures on 
$\cg$ is at least infinitesimally a
complex variety. 
Using the theory of 
deformations of complex structures, we are able to prove 
that for any small deformation of a rational complex structure $J$, 
the 
isomorphism \rf{iso} holds (Lemma 8).

\medskip

\noindent If $M$ is a compact complex parallelisable nilmanifold, i.e., $G$ is 
a nilpotent complex
Lie group and $J$ is also right invariant, we have that $\cg^{i}_J =  
\cg^i$ and Theorem B follows from \cite [Theorem 1] {Sak}.   

\noindent An
important class of complex structures is given by the abelian ones 
(i.e. those satisfying
the condition
$\lbrack J X, J Y \rbrack = \lbrack X, Y \rbrack$, for any $X, Y \in 
\cg$
\cite{BD,DF}). The nilmanifolds with
an abelian complex structure are to some extent  dual to complex
parallelisable nilmanifolds: indeed, 
 in the complex parallelisable case $$d \lambda^{1, 0}  \subset
\lambda^{2,0},$$ and in the abelian case $$ d \lambda^{1, 0}  \subset
\lambda^{1,1},$$  where $\lambda^{p, q}$ denotes the space of 
$(p,q)$-forms on
$\cg$. In this last case we will compute the minimal model of the 
Dolbeault
cohomology of $M$ and  prove  that the isomorphism \rf{iso} holds 
for any  abelian complex structure. In \cite{CFGU} was proved  a 
similar result for the Dolbeault cohomology of $M$ endowed with a 
nilpotent complex structure, which is a slight generalization of the 
abelian one. 

\smallskip

Note that, however, if $M$ is a complex solvmanifold $G/\Gamma$ (with $G$ a solvable not 
nilpotent Lie group), the isomorphism \rf{iso} does not hold 
in general, as shown in \cite{Le}; a discussion of the behaviour of 
the Dolbeult cohomology of homegeneous manifolds under group actions  
can be found in \cite{Ak}. 
If $G$ is a compact even dimensional (semisimple) Lie group  endowed with a left invariant complex 
structure, the Dolbeault cohomology of $G$ does not arise from just invariant 
classes, as the example in \cite{Pi} shows.

\smallskip

 This paper is organized as follows.

In Section \rl{Section 1}, following \cite{Sa},
we define a descending series  of subalgebras $\{ \cg_J^i \}$ (with 
$\cg_J^0 = \cg$ and $\cg_J^{s +
1} = \{ 0 \}$) for the Lie algebra $\cg$ associated to the complex 
structure $J$. In general, the
subalgebra $\cg_J^i$ is  not a rational subalgebra of $\cg^{i - 
1}_J$.  If  $J$ is rational,  any 
$\cg_J^i$  is rational in $\cg^{i - 1}_J$.

The importance of this series is twofold. 

\noindent First (Section \rl{Section 
2}), in the case that the subalgebra $
\cg_J^i$ is rational in $\cg^{i - 1}_J$ (in particular if $J$ is 
rational), it allows us  to define
a set of holomorphic fibrations of nilmanifolds:

\noindent $\tilde p_0: M = G/ \Gamma \to G^{0,1}_J/p_0 (\Gamma)$, 
with standard fibre
$G^1_J/\Gamma^1$,

$\ldots$

\noindent $\tilde p_{s - 1}: G^{s - 1}_J/ \Gamma^{s - 1} \to 
G^{s-1,s}_J/p_{s-1}
(\Gamma^{s - 1})$ with standard fibre $G^s_J/\Gamma^s$.

\noindent For the above fibrations we will consider the associated Borel 
spectral sequence $(E_r, d_r)$ 
(\cite[Appendix II by A. Borel, Theorem 2.1]{Hi}, \cite{FW,Lp}), 
which relates the Dolbeault
cohomology of each total space with the Dolbeault cohomology of each 
base and fibre.       
\smallskip

\noindent Secondly (Section \rl{Section 3}) following \cite{Sa}, we 
will prove that one can choose a
basis of $(1, 0)$-forms (and also $(0, 1)$-forms) on $\cg$ which is 
compatible
with the above descending series. This basis will give a basis of 
$G^{i - 1,
i}_J$-invariant $(1, 0)$-forms on the nilmanifolds $G^{i - 1, 
i}_J/p_{i-1}
(\Gamma^{i - 1})$ and of $G^i_J$-invariant $(1, 0)$-forms on the 
nilmanifolds
$G^i_J/\Gamma^i$, $i = 0, ..., s$.

\noindent Next, in Section \rl{Section 4}, we consider a spectral 
sequence $(\tilde E_r, \tilde d_r)$ 
concerning the $G^{i-
1}_J$-invariant Dolbeault cohomology  of each total space $G^{i -
1}_J/\Gamma^{i - 1}$, the $G^{i - 1, i}_J$-invariant Dolbeault 
cohomology  
each base $G^{i - 1, i}_J/p_{i - 1} (\Gamma^{i - 1})$ and  the
$G^i_J$-invariant Dolbeault cohomology  of each fibre 
$G^i_J/\Gamma^i$.
 In this way $(\tilde E_r, \tilde d_r)$ is relative to the Dolbeault 
cohomologies of the Lie
algebras
$\cg_J^i$,
$\cg_J^{i-1}$ and $\cg_J^{i-1}/\cg_J^i$. Note that the latter are the
underlying Lie algebras of the fibre, the total space and the base,
respectively, of the above holomorphic fibrations. 

 In Section \rl{Section 5} we compare the spectral sequence 
$(\tilde E_r, \tilde d_r)$ with the Borel
spectral sequence $(E_r, d_r)$.   
Inductively  (starting with $i = s$) these two spectral sequences  
allow us to
give isomorphisms between the Dolbeault cohomologies of the total 
spaces and
the one of the corresponding Lie algebras. The last step gives 
Theorem B.

\smallskip

Note that our construction of this set of holomorphic fibrations is in
the same vein as principal holomorphic torus towers, 
introduced by Barth and Otte
\cite{BO}. In some cases, like the one of abelian complex structures,
$G/\Gamma$ is really a principal holomorphic torus tower.

\smallskip

 In Section \rl{Section 6}  we will give a proof of Theorem A.

\smallskip

In Section \rl{Section finale} we give examples of compact nilmanifolds with non rational complex structures.

\medskip

\noindent We wish to thank Prof. Simon Salamon for the great deal of  
suggestions he
gave us and his constant encouragement. We are also grateful to Prof. 
Isabel
Dotti for useful conversations and the hospitality at the  FaMAF of 
C\' ordoba
(Argentina).  

\section{A descending series associated to the complex 
structure}\lb{Section 1}

 We recall that, since $G$ is $(s+1)$-step nilpotent, one has the
descending central series $\{ \cg^i \}_{i \geq 0}$, where $$ \cg =\cg^0
\supseteq \cg^1 =[\cg, \cg]\supseteq  \cg^2 =[\cg^1, \cg] \supseteq ...
\supseteq \cg^s \supseteq \cg^{s+1}=\{0\}.\eqno(D) $$ We define the following
subspaces of $\cg$  $$ \cg^i_J:=\cg^i+J \cg^i. $$ Note that $\cg^i_J$ is $J$
-invariant. \medskip

\noindent{\bf Lemma 1}. {\sl 

\noindent  (1) $\cg^i_J$ is an ideal of $\cg^{i-1}_J$.

\noindent (2) $\cg^{i-1}_J / \cg^i_J$ is an abelian algebra.

\noindent (3) $\cg^s_J$ is an abelian ideal of $\cg^{s-1}_J$.} \medskip

\noindent{\it Proof.}   
 
\noindent (1) For any $X = X_1 + J X_2 \in \cg^{i - 1}_J$ and $Y = Y_1 + J Y_2 \in
\cg^i_J$ (with $X_l \in \cg^{i - 1}$ and $Y_k \in \cg^i$), we have that $$
\lbrack X,  Y \rbrack = \lbrack X_1,  Y_1 \rbrack + \lbrack X_1,  J Y_2
\rbrack + \lbrack J X_2,  Y_1 \rbrack + \lbrack J X_2,  J Y_2 \rbrack. $$ We
can easily see  that $\lbrack X_1, Y_1 \rbrack$, $\lbrack X_1, J Y_2 \rbrack $
and $\lbrack J X_2, Y_1 \rbrack$  belong to $\cg^i$ by definition of the
descending central series. Moreover $\lbrack J X_2,  J Y_2 \rbrack$ belongs to
$\cg^i_J$ because $J$  satisfies an integrability condition, namely the
Nijenhuis tensor $N$ of $J$, given by $$ N (Z, W) =\lbrack Z, W \rbrack + J
\lbrack  J Z, W \rbrack + J \lbrack  Z,  J W \rbrack -  \lbrack  J Z,  J W
\rbrack, \quad Z, W \in \cg, $$ must be zero \cite{NN}.

\smallskip

\noindent (2) For any $X = X_1 + J X_2$ , $Y  = Y_1 + J Y_2$ elements of $\cg^{i -
1}_J$ we have that $$ \lbrack X + \cg^i_J, Y +  \cg^i_J \rbrack = \lbrack X_1
, Y_1 \rbrack  + \lbrack X_1, J Y_2 \rbrack + \lbrack  J X_2 , Y_1 \rbrack +
\lbrack J X_2 , J Y_2 \rbrack  + \cg^i_J. $$ Then using the same argument as
in  (1)  it follows that $\lbrack X + \cg^i_J, Y +  \cg^i_J \rbrack =
\cg^i_J.$ 

\smallskip

\noindent (3) Using the fact that $\cg^s$ is central (i.e. $\lbrack \cg^s, \cg
\rbrack = 0$) and that $N = 0$ it is possible to prove that $\lbrack X,  Y
\rbrack$ vanishes for any $X, Y \in \cg^s_J$. \qed

\medskip

\noindent Observe moreover that any $\cg^i_J$ is nilpotent. 

\smallskip

Hence we have the descending series $$ \cg=\cg^0_J \supset \cg^1_J
\supseteq  \cg^2_J \supseteq ... \supseteq \cg^s_{J} \supseteq
\cg^{s+1}_J=\{0\}.\eqno(DJ) $$

\noindent {\it Remark 1.} The first inclusion $\cg^1_J \subset \cg$ is always strict
\cite [Corollary 1.4] {Sa}. 

\noindent Observe also that in case of complex parallelisable nilmanifolds 
\cite{Sak}, the
filtration $\{ \cg^i_j\}$ coincides with the descending central series $\{
\cg^i \}$ and then the $\{ \cg^i_j\}$ are rational.

\noindent In general, given a rational structure
${\cg}_{\Q}$ for $\frak g$,  we say that a $\R$-subspace $\frak h$   of $\cg$ is rational
if ${\ch}$ is the $\R$-span of ${\ch}_{\Q} = {\ch} \cap {\cg}_{\Q}$.
In general $\cg^1_J$ is not a rational subalgebra of $\frak g$.
 When $J$ is rational, it is possible to prove  that  $\cg^i_J$ 
is rational in $\cg^{i - 1}_J$. Indeed,  we have that $\cg^i$ is 
rational in $\cg^{i-1}$. Then 
 $\cg^i = \R{\rm {-span}} \{\cg^i \cap \cg^{i - 1}_{\Q} \}$. Since 
$J \cg^{i - 1}_{{\Q}} \subseteq \cg^{i - 1}_{{\Q}}$ it follows that $\cg^i_J = \R{\rm {-span}} \{\cg^i_J \cap 
\cg^{i - 1}_{\Q}\}$. 
 Moreover when $J$ is abelian,
$\cg^i_J$ is an ideal of
$\cg$, for any $i$ and the center
$$ \cg_1 = \{ X \in \cg \mid \lbrack X, \cg \rbrack = 0 \} $$
is a rational $J$-invariant ideal of $\cg$.

\section{Holomorphic fibrations and Borel spectral 
sequences}\lb{Section 2}

In this section we suppose that the complex structure  $J$ is rational 
and we associate  a set of holomorphic
fibrations to the above  descending series. We recall that a {\it holomorphic fibre bundle} $\pi: T
\rightarrow B$ is a a holomorphic map between the complex manifolds $T$ and
$B$,  which is locally trivial,   
whose typical fibre $F$ is a complex manifold 
and such that the transition functions are holomorphic. By definition the structure 
group  (i.e. the group of holomorphic automorphisms of the typical fibre) is
a complex Lie group.

 To define the above fibrations, we consider first the surjective homomorphism
$$  p_{i - 1}: \cg^{i - 1}_J \rightarrow \cg^{i - 1}_J  / \cg^i_J, $$  for
each $i = 1, \ldots, s$. If  $G^{i - 1}_J$ and $G^{i-1, i}_J$ denote
 the simply connected nilpotent Lie group corresponding to $\cg^{i - 1}_J$
 and $\cg^{i - 1}_J  / \cg^i_J$ respectively, we have the surjective
homomorphism $$ p_{i - 1}: G^{i - 1}_J \rightarrow G^{i-1, i}_J. $$

\noindent We define inductively $G^i_J$ to be the fibre of $p_{i - 1}$. Remark that
the Lie algebra of $G^i_J$ is $\cg^i_J$.

Given the uniform discrete subgroup $\Gamma$ of $G=G^0$,  we consider the
continuous surjective map $$ \tilde p_0: G/ \Gamma \rightarrow G^{0, 1}_J /
p_0 (\Gamma). $$   Since $J$ is rational, $\cg^1_J$ is a rational subalgebra of $\frak g$, then
$\Gamma^1:=\Gamma
\cap G^1_J$ is a uniform discrete subgroup of $G^1_J$ \cite[Theorem 
5.1.11]{CG}.  Then, by \cite[Lemma
5.1.4 (a)]{CG}, 
$p_0 (\Gamma)$ is a  a uniform discrete subgroup  of $G^{0, 1}_J$ (i.e.  $G^{0, 1}_J /
p_0 (\Gamma)$ is compact, cf. \cite{Ra}).  

\noindent Note moreover that $G^1_J$ is simply connected. This follows from
the homotopy exact sequence of the fibering $p_0$. Indeed we have $$ \ldots
\rightarrow \pi_2(G^{0, 1}_J)=(e) \rightarrow \pi_1 (G^1_J)  \rightarrow \pi_1
(G)=(e) \rightarrow \ldots $$ Finally it is not difficult to see that $G^1_J$
is connected. Indeed, if $C$ is the connected component of the identity in
$G^1_J$,  id$: G \to G$ induces a  covering homomorphism $G/C \to G/G^1_J\cong
G^{0,1}_J$ which must be the identity, since $G^{0,1}_J \cong \R^{N_0}$. Thus
$C=G^1_J$.

\medskip

Now one can  repeat the same construction for any $i$, since $\cg^i_J$ is a rational ideal of
$\cg^{i - 1}_J$. So, for any
$i=1,
\ldots , s$ we have a
 map $$ \tilde p_{i-1}: G^{i-1}_J / \Gamma^{i-1} \to G^{i-1, i}_J / p_{i-1}
(\Gamma^{i-1}). $$

\medskip

\noindent{\bf Lemma 2}. {\sl $\tilde p_{i-1}: G^{i-1}_J / \Gamma^{i-1} \to G^{i-1,
i}_J / p_{i-1} (\Gamma^{i-1})$ is a holomorphic fibre bundle. } \medskip

\noindent{\it Proof.} Observe first that $\tilde p_{i-1}$ is the induced map of
$p_{i-1}$ taking quotients of discrete subgroups.  The tangent map of $\tilde 
p_{i - 1}$  $$ \cg^{i-1}_J \rightarrow \cg^{i-1}_J / \cg^{i}_J $$ is
$J$-invariant. Thus $\tilde  p_{i - 1}$ is a holomorphic submersion. In
particular it is a holomorphic family of compact complex manifolds in the
terminology of \cite{KS} (see also \cite{Sun}). The fibres of $\tilde p_{i - 1}$ are all
holomorphically equivalent to $G_J^i/ \Gamma^i$ (the typical fibre). Thus a
theorem of Grauert and Fisher \cite{FG} applies, implying that $\tilde p_{i-1}$ is a
holomorphic fibre bundle.  \qed

\medskip

\noindent Note that $G^{i-1}_J / \Gamma^{i-1}$, $G^i_J/\Gamma^i$, $G^{i-1, i}_J /
p_{i-1} (\Gamma^{i-1})$ are compact connected nilmanifolds.

\medskip

Given a holomorphic fibre bundle it is possible to construct the associated
Borel spectral sequence, that relates the Dolbeault cohomology of the total
space $T$ with that of the basis $B$ and of the fibre $F$.  We will need the
following Theorem (which follows from \cite[Appendix II by A. Borel, Theorem
2.1]{Hi} and \cite{FW}).

\medskip

\noindent{\bf Theorem 3}.  {\sl Let $p:T \to B$ be a holomorphic fibre bundle, with
compact connected fibre $F$  and $T$ and $B$ connected. Assume that either

\noindent (I) $F$ is K\" ahler 

or

\noindent (I') the scalar cohomology bundle  $$ {\bf H}^{u, v} (F) = 
\bigcup_{b \in B}
\HD{u}{v} (p^{-1} (b))$$ is trivial.

Then there exists a spectral sequence $(E_r, d_r)$, ($r\geq 0$) with the
following properties:

\noindent (i) $E_r$ is 4-graded by the fibre degree, the base degree and the type.
Let $^{p,q}E_r^{u,v}$ be the subspace of elements of $E_r$ of type $(p,q)$,
fibre degree $u$ and base degree $v$. We have $^{p,q}E_r^{u,v}=0$ if $p+q\neq
u+v$ or if one of $p,q,u, v$ is negative. The differential $d_r$ maps 
$^{p,q}E_r^{u,v}$ into $^{p,q+1}E_r^{u+r,v-r+1}$.

\noindent (ii) If $p+q=u+v$  $$ ^{p,q}E_2^{u,v} \cong \sum_k    \HD{k}{u-k} (B) 
\otimes \HD{p-k}{q-u+k} (F). $$
 
\noindent (iii) The spectral sequence converges to $\hd (T)$.}

\section{An adapted basis of $(1, 0)$-forms} \lb{Section 3}

In this section we prove that one can choose a basis of
(1,0)-forms on $\cg$ which is compatible with the descending series $(DJ)$. We
consider, like in \cite{Sa}, some subspaces $V_i$ ($i=0, \ldots , s+1$) of $V:=(T_e
G)^*\cong \cg^*$, that determine a series, which is related to the descending
central series $(D)$.

\noindent Indeed we define: 
\[
\ba{lll} 
&V_0=\{0\}\cr &V_1=\{ \alpha \in V \mid
d\alpha =0\} \\
 &\dots\cr &V_i=\{ \alpha \in V \mid d\alpha \in \Lambda^2
V_{i-1} \} \\ &\dots \\ &V_{s+1}=V. \ea\]
Note that $V_i$ is the
annihilator $(\cg^i)^o$ of the subspace $\cg^i$ and that $\{0\}=V_0 \subseteq
V_1 \subseteq \ldots \subseteq V_{s+1}=V$ \cite[Lemma 1.1]{Sa}.

\noindent If we now let $(\cg^i_J)^0 \cap \lambda^{1,0}=: V_i^{1,0}$, by 
\cite[Lemma
1.2]{Sa} we have that there exists a basis of (1,0)-forms $\{ \om_1, \ldots \om_n\}$
such that if $\om_l \in V_i^{1,0}$ then $d \om_l$ belongs to the ideal (in
$(\cg^\C)^*$) generated by $V_{i-1}^{1,0}$ \cite[Theorem 1.3]{Sa}. In particular
there exists at least a closed $(1, 0)$-form (this implies Remark 1).

\noindent Moreover we have the following isomorphisms: $$ \left ( \cg_J^{i-1} / 
\cg_J^{i} \right )^\C \cong V_i^{1,0}/V_{i-1}^{1,0} \oplus
V_i^{0,1}/V_{i-1}^{0,1}, \qquad i=1,\ldots , s, $$ where $V_i^{0,1}$ is the
conjugate of $V_i^{1,0}$.

\noindent With respect to the subspaces $V_i^{1,0}$ the above basis can be ordered as
follows (we let $n_i:=\dim_\C \cg^i_J$):

\noindent $\om_1, \ldots , \om_{n-n_1}$ are elements of $V_1^{1,0}$  (such that  $d
\om_l = 0$) or $\cg/\cg_J^1$ is the real vector space underlying $V^{1,0}_1$;

\noindent $\om_{n-n_1+1}, \ldots , \om_{n-n_2}$ are elements of  $V_2^{1,0} \backslash
V_1^{1,0}$ or $\cg_J^1/\cg_J^2$ is the real vector space underlying the
quotient $V^{1,0}_2/V^{1,0}_1$;

$\ldots$

\noindent $\om_{n-n_{s-1}+1}, \ldots , \om_{n-n_s}$ are elements of  $V_s^{1,0}
\backslash V_{s-1}^{1,0}$;

\noindent $\om_{n-n_s+1}, \ldots , \om_{n}$ are elements of  $\lambda^{1,0} \backslash
V_s^{1,0}$. 

\smallskip

\noindent Here $V_i^{1,0} \backslash V_{i-1}^{1,0}$ denotes a complement of
$V_{i-1}^{1,0}$ in  $V_i^{1,0}$ (which corresponds to the choice of a
complement of $\cg^i_J$ in $\cg^{i-1}_J$).

\noindent Hence, by definition, the elements of $V_i^{1,0} \backslash V_{i-1}^{1,0}$
and, by identification, the elements of the quotient $V_i^{1,0} /
V_{i-1}^{1,0}$, are (1,0)-forms on $\cg$ which vanish on $\cg^i_J$. So they
may be identified with forms on the quotient $\cg^{i-1}_J/ \cg^i_J$.

\noindent In this way we can consider: 
 
\noindent - the elements of  $\lambda^{1,0} / V_1^{1,0} = V_2^{1,0} / V_1^{1, 0} 
\oplus V_3^{1,0} / V_2^{1, 0} \oplus \ldots \oplus \lambda_{1,0} / V_s^{1, 0}$
as  (1,0)-forms on $\cg^1_J$, 

 ...

\noindent - the elements of  $\lambda^{1,0} / V_{s - 1}^{1,0} =  V_s^{1,0} / V_{s -
1}^{1,0} \oplus \lambda^{1,0}/ V_s^{1, 0}$ as 
 (1,0)-forms on  $\cg^{s - 1}_J$ and 

\noindent - the elements of $\lambda^{1, 0}/ V_s^{1,0}$ as (1,0)-forms on $\cg^s_J$.

\smallskip

\noindent Thus we can prove  a Lemma on the existence of a basis of $(1,0)$ forms on
$\cg$ related to the series $(DJ)$.

\medskip

\noindent{\bf Lemma 4}.  {\sl It is possible to choose a basis of (1,0)-forms   on
$\cg$ such that (with respect to the order of before)  $\{ \om_{n - n_{i - 1} +
1}, \ldots \om_{n - n_i}, \ldots \om_n \}$ is a basis of (1,0)-forms on $\cg^{i
- 1}_J$. Moreover we can consider (up to identifications)
 $\{\om_{n - n_i + 1}, \ldots \om_n \}$ as forms on $\cg^i_J$ and 
 $\{\om_{n - n_{i - 1} + 1}, \ldots \om_{n - n_i} \}$ as forms on   $\cg^{i-1}_J
/ \cg^i_J$. } 

\medskip

\noindent{\it Proof.} By the above arguments, for any $i =1, \ldots,  s+1$, it is
possible to choose a basis of (1,0)-forms on $\cg^{i - 1}_J$ as elements of
$\lambda^{1, 0}/ V_{i - 1}^{1, 0} = \lambda^{1, 0}/ V_i^{1, 0} \oplus V_i^{1,
0} / V_{i - 1}^{1, 0}$.  With respect to the above decomposition  the forms
on  $\lambda^{1, 0}/ V_i^{1, 0}$ can be identified with forms on $\cg^i_J$ 
extended by zero on $\cg^{i - 1}_J$ and the forms on $V_i^{1, 0} / V_{i -
1}^{1, 0}$  with forms on $\cg^{i - 1}_J / \cg^i_J$, because these forms
vanish on $\cg^i_J$. \qed

\medskip

\noindent{\it Remark 2.} $d\om_i$, $i=n - n_{i - 1} + 1, \ldots , n - n_i$ belongs to
the ideal generated by $\{ \om_l, \quad l=1, \ldots ,n - n_{i - 1} \}$.
 
\noindent{\it Remark 3.} If $J$ is abelian it is possible to choose a basis of $(1,
0)$-forms 
 $\{ \om_1, \ldots, \om_n \}$ on  $\cg$ such that  $$ d\om_i \in 
 \wedge^2 \left< \om_1,
\ldots, \om_{i - 1}, \overline \om_1, \ldots, \overline \om_{i - 
1}\right> \cap \,
\lambda^{1, 1}.$$ 
 
\section {A spectral sequence for the complex of invariant forms} 
\lb{Section 4}
 
We construct a spectral sequence $^{p,q}\tilde E_r^{u,v}$ for the
complexes of $G^{s-1}_J$-invariant forms on $G^{s-1}_J / \Gamma^{s-1}$ whose
Dolbeault cohomology identifies with $\hd ((\cg^{s-1}_J)^\C)$.  To do
this, we give a filtration of the complex $\Lambda_{\ct} = \oplus \Lambda_{\ct}^{p,q}$ of 
differential forms of type $(p,q)$ on  $\ct = {\cg}^{i-1}_J$. 

\noindent We know from Lemma 4 that there exists a basis of $(1,0)$-forms
$\om_h^{\ct}$ on $\ct$ (and of $(0,1)$-forms $\overline \om_h^{\ct}$) such that
part of them are  $(1,0)$-forms $\omega_j^{\cb}$ on ${\cb} = \cg^{i -
1}_J/{\cg}^i_J$ and part are  forms $\om_k^{\cf}$ on $\cf = \cg^i_J$. We define 
\[
\ba{ll}
\tilde L_k:=&\{ \om^{\ct} \in \Lambda_{\ct} \mid \om^{\ct} {\hbox{
is a sum of monomials }}\\ & \om_I^{\cb} \w \overline \om_J^{\cb} \w
\om_{I'}^{\cf} \w \overline \om_{J'}^{\cf} {\hbox { in which }} |I|+|J| \geq k
\}, \ea
\]
 where $|A|$ denotes the number of elements of the finite set $A$.

\noindent Note that $\tilde L_0= \Lambda_{\ct}$ and that 
\[
\ba{ll}&\tilde L_k=0 \qquad {\hbox{for}} \quad k > \dim_\R \cb,\\
 &\tilde
L_k \supset \tilde L_{k+1}, \qquad \debar \tilde L_k \subseteq \tilde L_k,
\quad k \geq 0 \ea
\] The above shows that $\{\tilde L_k \}$ defines a
bounded decreasing filtration of the differential module $(\Lambda_{\ct},
\debar)$. 

\noindent Of course  $$ \tilde L_k=\sum_{p,q} {^{p,q}{\tilde L}_k},  {\hbox { where
}} \quad ^{p,q}{\tilde L}_k=\tilde L_k \cap \Lambda_{\ct}^{p,q}  $$ and the
filtration is compatible with the bigrading provided by the type (and also
with the total degree).

Recall that, by definition (see e.g. \cite{GH}) $$ ^{p,q}{\tilde E}^{u,v}_r={{
^{p,q} Z^{u, v}_r } \over { ^{p,q}Z^{u+1,v-1}_{r-1} +  {^{p,q}B}^{u,v}_{r-1}
}}, $$ where  
\[
\ba{ll} & ^{p,q}Z^{u,v}_r={^{p,q}\tilde L}_u
(\Lambda^{u+v}_{\ct}) \cap \ker \debar ({^{p,q+1}\tilde L}_{u+r}
(\Lambda^{u+v+1}_{\ct})) \\ & ^{p,q}B^{u,v}_r={^{p,q}\tilde L}_u
(\Lambda^{u+v}_{\ct}) \cap
 \debar ({^{p,q-1}\tilde L}_{u+r} (\Lambda^{u+v-1}_{\ct})) \ea
 \]
Moreover (cf. \cite{GH}) $$ ^{p,q}\tilde E_0^u={{ ^{p,q} \tilde L_u } \over {
^{p,q}\tilde L_{u+1} }}, $$ where we denote by $^{p,q} \tilde E^{u,v}_r$ and
$^{p,q}\tilde E^{u}_r$ the spaces of elements of type $(p,q)$ and total degree
$u+v$ and degree $u$ respectively in the grading defined by the filtration.

\noindent Note also that an element of $\tilde L^k$ identifies with an element of
$\sum_{c + d \geq k} \Lambda^{a,b}_{\cf} \otimes \Lambda^{c,d}_{\cb}$.

\medskip

\noindent{\bf Lemma 5}.  {\sl Given  the holomorphic fibration $$G^{i - 1}_J/\Gamma^{i - 1}
\rightarrow G^{i-1, i}_J /p_{i - 1} (\Gamma^{i - 1}),$$  with standard fibre $G^i_J/\Gamma^i$,
the spectral sequence $(\tilde E_r, \tilde d_r)$ ($r \geq 0$) converges to
$\hd ((\cg^{i-1}_J)^\C)$ and
 $$ ^{p,q}\tilde E_2^{u,v} \cong \sum_k \HD{k}{u-k} ((\cg^{i-1}_J/\cg^i_J)^\C)
\otimes \HD{p-k}{q-u+k} ((\cg^i_J)^\C). \eqno (\tilde Ii)$$}

\noindent{\it Proof.}
The fact that $(\tilde E_r, \tilde d_r)$ converges to $\hd ((\cg^{i-1}_J)^\C)$ is a
general property of spectral sequences associated to filtered complexes (cf.
\cite{GH}).

Let $[\om ] \in \,  ^{p,q}{\tilde E}_0^u={{^{p,q} \tilde L_u} \over {^{p,q}\tilde L_{u+1}}}$.  We
compute the differential $\tilde d_0: {^{p,q}\tilde E}^u_0 \rightarrow  {^{p,q+1}\tilde E}^u_0$
defined by
$\tilde d_0 [\om]= [\debar \om]$.   We can write (up the above identifications) $\om=\sum
\om_I^{\cb} \w \overline \om_J^{\cb} \w \om_{I'}^{\cf} \w \overline
\om_{J'}^{\cf}$ where $|I|+|J|=u$  (because we operate $\mod \tilde L_{u+1}$),
 and $|I'|+|J'|=p+q-u$, since  $|I'|+|J'|+|I|+|J|=p+q$. Moreover,  using the
fact that $\debar$ sends forms in ${\cb}$ on forms that either are in $\cb$ or
vanish on $\cb$ (cf. Remark 2) and that $\cb$ is abelian, we get  
\[
\ba{ll}  \debar \om=&\sum
{\debar}_{\cb} (\om_I^{\cb} \w \overline \om_J^{\cb}) \w (\om_{I'}^{\cf} \w
\overline \om_{J'}^{\cf}) + (-1)^s  (\om_I^{\cb} \w \overline \om_J^{\cb}) \w
\debar (\om_{I'}^{\cf} \w \overline \om_{J'}^{\cf})=\\
 =& (-1)^s 
(\om_I^{\cb} \w \overline \om_J^{\cb}) \w {\debar}_{\cf} (\om_{I'}^{\cf} \w
\overline \om_{J'}^{\cf}) \mod \tilde L_{u+1},\ea
\]
where ${\debar}_{\cb}$ and ${\debar}_{\cf}$ denote the differential on the complexes 
$\Lambda_{\cb}$ and  $\Lambda_{\cf}$, respectively.  Thus  $$\tilde d_0
[\om] = [\debar_{\cf} \om],$$ which implies  $$
^{p,q}\tilde E^u_1 \cong  \sum_k   \HD {p-k}{q-u+k} (\cf^{\C})
\otimes \Lambda^{k,u-k}_{\cb}. \eqno (k1) $$ 

\noindent Performing the same proof as in \cite[Appendix II by Borel, Section 6]{Hi} and using the fact that
$\cb$ is abelian, it is possible to prove that
$\tilde d_1$ identifies with
$\debar$ via the above isomorphism and that we have  $$ ^{p,q}\tilde
E^u_2 \cong  \sum_k  \HD {k}{u-k} (\cb^{\C}) \otimes \HD {p-k}{q-u+k} (\cf^{\C})
 ( \cong \sum_k  \Lambda_{\cb}^{k,u-k}  \otimes \HD {p-k}{q-u+k} (\cf^{\C})).  \eqno (k2) $$ 
\qed

\bigskip

\section {Proof of Theorem B} \lb{Section 5}

\noindent First we note that Theorem B is trivially true if the
nilmanifold comes from an abelian group, i.e., it is a complex torus. Namely,
if $A/\Gamma$ is a complex torus we have $$ \hd(A/\Gamma)\cong \hd (\ca^\C).
\eqno (a) $$

\noindent We consider the holomorphic fibrations 
\[
\ba{ll}
&G^s_J/\Gamma^s \hookrightarrow G^{s - 1}_J/\Gamma^{s - 1} \rightarrow G^{s -
1, s}_J/p_{s - 1} (\Gamma^{s-1})\cr &\ldots \\
 &G^1_J/\Gamma^1
\hookrightarrow G/\Gamma  \rightarrow G^{0, 1}_J/p_0 (\Gamma).
\ea
\]
 The aim
is to obtain  informations about the Dolbeault cohomology of $G/\Gamma$
inductively through the Dolbeault cohomologies of $G^i_J/\Gamma^i$ (the
nilmanifolds $G^i_J/\Gamma^i$ play alternately the r\^ oles of fibres and
total spaces of the above fibre bundles).  Note that since the bases are
complex tori, $\hd (G^{i, i - 1}_J/p_{i - 1} (\Gamma^{i - 1}) \cong \hd ( (
\cg^{i - 1}/\cg^i)^{\C}).$  
 
 To this purpose we will associate to these
fibrations two spectral sequences. 

\noindent The first is a version of the Borel spectral sequence (considered 
in Section \rl{Section 2}) which relates
the Dolbeault cohomologies of the total spaces with those of fibres and bases. 

\noindent The second is the spectral sequence $(\tilde E_r, \tilde d_r)$ (constructed in the
previous Section) relative to the Dolbeault cohomologies of the Lie algebras
$\cg^i, \cg^{i - 1}, \cg^{i - 1}/\cg^i$.

\smallskip

We will proceed inductively on the index $i$ in the descending series
$(DJ)$, starting from $i=s$.
 
\smallskip 

\noindent {\it First inductive step.} Let us use the holomorphic fibre bundle
$$\tilde p_{s-1}: G^{s-1}_J / \Gamma^{s-1} \to G^{s-1, s}_J / p_{s-1}
(\Gamma^{s-1})$$ with typical fibre $G^{s}_J / \Gamma^{s}$. 

\noindent Recall (cf. Lemma 1) that $G^{s}_J / \Gamma^{s}$ and $G^{s-1, s}_J /
p_{s-1} (\Gamma^{s-1})$ are complex tori. Thus by $(a)$ 
\[
\ba{lll}
& \hd
(G^{s}_J / \Gamma^{s})\cong \hd ((\cg^s_J)^\C) &(s) \\
 & \hd (G^{s-1, s}_J /
p_{s-1} (\Gamma^{s-1}))\cong \hd ((\cg^{s-1}_J/\cg^s_J)^\C). & (s,s-1)
\ea
\]

\medskip

\noindent Applying Theorem 3  (since the fibre $G^{s}_J / \Gamma^{s}$ is K\" ahler) 
and using $(s)$ and $(s,s-1)$, we get $$ ^{p,q}E_2^{u,v} \cong \sum_k
\HD{k}{u-k} ((\cg^{s-1}_J/\cg^s_J)^\C) \otimes \HD{p-k}{q-u+k} ((\cg^s_J)^\C).
\eqno (Is)$$

\smallskip

Next we use the spectral sequence $(\tilde E_r, \tilde d_r)$. Note that the inclusion between the
Dolbeault complex of $G^{s-1}_J$-invariant forms on $G^{s-1}_J / \Gamma^{s-1}$ and the forms on
$G^{s-1}_J /
\Gamma^{s-1}$, induces an inclusion of each term in the spectral sequences $$
^{p,q}\tilde E_r^{u,v} \subseteq ^{p,q}{E}_r^{u,v}, $$ (which is actually a
morphism of spectral sequences).

\noindent By Lemma 5, for $i=s$, we have that $(\tilde E_r, \tilde d_r)$ converges to $\hd
((\cg^{s-1}_J)^\C)$ and 
 $$ ^{p,q}\tilde E_2^{u,v} \cong \sum_k \HD{k}{u-k} ((\cg^{s-1}_J/\cg^s_J)^\C)
\otimes \HD{p-k}{q-u+k} ((\cg^s_J)^\C). \eqno (\tilde Is)$$

Comparing $(Is)$ with $(\tilde Is)$, we get that $E_2 = \tilde
E_2$, hence the spectral sequences $(E_r, d_r)$ and $(\tilde E_r, \tilde d_r)$ converge to the same
cohomologies. Thus
 $$
 \hd (G^{s-1}_J / \Gamma^{s-1})\cong \hd ((\cg^{s-1}_J)^\C). \eqno (s-1)  $$

\medskip

\noindent{\it General inductive step.} We use the holomorphic fibre bundle $$\tilde
p_{i-1}: G^{i-1}_J / \Gamma^{i-1} \to G^{i-1, i}_J / p_{i-1} (\Gamma^{i-1})$$
with typical fibre $G^{i}_J / \Gamma^{i}$.  We assume inductively that
 $$ \hd (G^{i}_J / \Gamma^{i})\cong \hd ((\cg^i_J)^\C). \eqno (i) $$

\medskip

\noindent{\bf Lemma 6}. {\sl The scalar cohomology bundle  $$ {\bf H}^{u, v}
(G^i_J/\Gamma^i)  = \bigcup_{b \in G^{i-1,i}_J /p_{i - 1} (\Gamma^i)} \HD{u}{v}
({\tilde p}_{i - 1}^{-1} (b))$$ is trivial.}

\medskip

\noindent{\it Proof.} By \cite[Section 5, formula 5.3]{KS} there exists a locally finite
covering  $\{ U_l \}$  of  $G^{i-1,i}_J/p_{i - 1} (\Gamma^{i - 1})$ such that
the action of the structure group of the holomorphic fiber bundle on $U_l \cap
(G^i_J/\Gamma^i)$  is the differential of the change of complex coordinates on
the fibre, so one can restrict oneself to consider the left translation by
elements of $G^{i - 1}_J$ as change of coordinates. Then the  scalar cohomology
bundle  ${\bf H}^{u, v} (G^i_J/\Gamma^i)$ is trivial  since  any of its 
fibres  is canonically isomorphic to  $\hd ((\cg^i_J)^{\C})$. More
explicitly,  a global frame for ${\bf H}^{u, v} (G^i_J/\Gamma^i)$ is given as
follows: for any cohomology class  $$ \alpha \in \HD{u}{v} ({\tilde p}_{i -
1}^{-1}  (1)) =  \HD{u}{v} (G^i_J/\Gamma^i) \cong \HD{u}{v}
((\cg^i_J)^{\C}),$$ ($1$: identity element of $G^{i-1,i}_J/p_{i - 1}
(\Gamma^i)$) one can take the corresponding cohomology class $\omega \in \HD{u}{v}
((\cg^i_J)^{\C})$ and regard it as a $G^{i - 1}_J$-invariant differential form
on $G^{i - 1}_J/\Gamma^{i - 1}$. Thus $b \mapsto \omega \vert_{{\tilde p}_{i -
1}^{-1} (b)}$ gives a global holomorphic  section of ${\bf H}^{u, v}
(G^i_J/\Gamma^i)$.  Taking a basis of $\hd ((\cg^i_J)^{\C})$ one gets a
global holomorphic frame of ${\bf H}^{u, v} (G^i_J/\Gamma^i)$. \qed

\medskip

Thus the assumption (I') in Theorem 3 is fulfilled. 

\noindent Observe moreover that, since $G^{i-1, i}_J / p_{i-1} (\Gamma^{i-1})$ is a
complex torus, by $(a)$,  $$ \hd (G^{i-1, i}_J / p_{i-1} (\Gamma^{i-1}))\cong
\hd ((\cg^{i-1}_J/\cg^i_J)^\C). \eqno (i,i-1)$$ Hence, by Theorem 3, we have $$
^{p,q}E_2^{u, v} \cong \sum_k \HD{k}{u-k} ((\cg^{i-1}_J/\cg^i_J)^\C) \otimes
\HD{p-k}{q-u+k} ((\cg^i_J)^\C). \eqno (Ii) $$ Then we proceed exactly like in
the first inductive step. Namely, we consider the spectral sequence
$^{p,q}\tilde E_r^{u,v}$ for the complexes of $G^{i-1}_J$-invariant forms on
$G^{i-1}_J / \Gamma^{i-1}$ whose Dolbeault cohomology identifies with $\hd
((\cg^{i-1}_J)^\C)$.  By Lemma 5 we have that $(\tilde E_r, \tilde d_r)$ converges to $\hd
((\cg^{i-1}_J)^\C)$ and
 $$ ^{p,q}\tilde E_2^{u,v} \cong \sum_k \HD{k}{u-k} ((\cg^{i-1}_J/\cg^i_J)^\C)
\otimes \HD{p-k}{q-u+k} ((\cg^i_J)^\C). \eqno (\tilde Ii)$$

\medskip

\noindent Using $(Ii)$ and proceeding like in the first inductive step we get the
proof of Theorem B. \qed

\bigskip

\noindent {\bf Remark on abelian complex structures. }
If the invariant structure $J$ is abelian, the centre $\cg_1$ is a rational $J$-invariant ideal
of $\cg$ (cf. Remark 1 in Section \rl{Section 1}). 

We give an alternative (and simpler) proof of that given in \cite{CFG} in the abelian case.  In the
same vein of \cite{No} we consider the principal holomorphic fibre bundle
$$
G/\Gamma \to G/ (\Gamma G_1),
$$ 
with typical fibre $G_1/(\Gamma \cap G_1) \cong \Gamma G_1/ G_1$ (where $G_1$ is the simply
connected Lie group corresponding to $\cg_1$).

\noindent First we can consider the Borel spectral sequence $(E_r, d_r)$ associated to the
Dolbeault complex $\Lambda^{*,*} (G/\Gamma)$ (cf. Theorem 3) and
a  spectral sequence  $(\tilde E_r, \tilde d_r)$ associated to $\Lambda^{*,*} (\cg)^\C$ and
constructed like in Section \rl{Section 4}. As to the latter spectral sequence $(\tilde E_r, \tilde d_r)$,
observe that it is possible to prove that there exists a basis 
$\{ \om_1, \ldots , \om_n \}$ of $(1,0)$- forms
on $\cg$  such that $\{ \om_1, \ldots , \om_{n-k}\}$
is a basis on $\cg/\cg_1$ (with $\dim \cg_1 = 2 k$) and $\{ \om_{n-k+1}, \ldots , \om_{n}\}$ is a
basis on $\cg_1$. So one can perform the same construction as in Section 
\rl{Section 4} (with $\ct = \cg$, $\cb
= \cg / \cg_1$ and $\cf = \cg_1$), using the fact that
$\cf$ is abelian. Thus one gets 
$$ ^{p,q}\tilde E_2^{u,v} \cong \sum_k \HD{k}{u-k} ((\cg/\cg_1)^\C)
\otimes \HD{p-k}{q-u+k} ((\cg_1)^\C)$$
and $(\tilde E_r, d_r)$ converges to $\hd ((\cg)^\C)$. 

\noindent Moreover, 
since $G_1/\Gamma \cap G_1$ and $G/(\Gamma G_1)$ are complex  torus, $\hd (G_1/\Gamma \cap
G_1) \cong \hd ((\cg_1)^\C)$ and $\hd (G/(\Gamma G_1) \cong \hd ((\cg/\cg_1)^\C)$ , so  by Theorem 3
we have $$ ^{p,q} E_2^{u,v} \cong \sum_k \HD{k}{u-k} (G/(\Gamma G_1))
\otimes \HD{p-k}{q-u+k} ((\cg_1)^\C).$$
This implies that 
Theorem B holds also for complex abelian structures. 

\medskip

Next we  construct  a minimal model for the Dolbeault cohomology of a
nilmanifold endowed with an abelian complex structure.

Recall that a {\it model} for the Dolbeault cohomology of $M$ is a differential
bigraded algebra $(\MM^{*,*}, \debar)$ for which there
exists a homomorphism $\rho: \MM^{*,*} \to \Lambda^{*,*} (M)$ of differential 
bigraded algebras inducing an isomorphism $\rho^*$ on the respective 
Dolbeault cohomologies  \cite{NT}.

\noindent Suppose $\MM$ is free on a vector space $V$. Then $\debar$ is called
decomposable if there is an ordered basis of $V$ such that the differential
$\debar$ of any generator $v$  of $V$ can be expressed in terms of the
elements of the basis preceding $v$. A model is called {\it minimal} if
$\MM$ is free and $\debar$ is decomposable \cite{Su}.

By \cite{Sa}, in the abelian case, there exists a basis $\{ \om_1, \ldots, \om_n \}$
of $(1, 0)$-forms on $\cg$
 such that
$$
d \om_i = \sum_{j <  k < i} A_{ijk}\,  \om_j  \wedge \overline \om_k, \quad i =
1, \ldots, n. \eqno{(m)} $$
Thus, by  Theorem B and $(m)$,   $(\Lambda^{*, *} (\cg^{\C}), \debar)$
is
 a minimal model for the Dolbeault cohomology of
$G/\Gamma$.

\section {Proof of Theorem A} \lb{Section 6}

Let
$$
{\cal C} ({\cg}) = \{ J \in End (\cg) \mid J^2 = - id, \lbrack J X, J Y \rbrack = \lbrack X, Y
\rbrack + J \lbrack JX, Y \rbrack + J \lbrack X, J Y \rbrack \}
$$
denote the set of complex structures on $\cg$. We will use the same 
notation as in \cite{Sa}. If $M = G/ \Gamma$ is a nilmanifold associated to
$\cg$, then
$$
0 \to Hom (\lambda^{1,0}, \lambda^{0,1}) \stackrel{\debar} \to Hom (\lambda^{1,0}, \lambda^{0,2})
 \stackrel{\debar} \to \ldots \stackrel{\debar} \to Hom (\lambda^{1,0}, \lambda^{0,n}) \to 0
$$
is a subcomplex of the Dolbeault complex of $M$ tensored with the holomorphic tangent bundle
$T^{1,0} M$. The Kernel $K$ of $\debar$ acting on $Hom (\lambda^{1,0}, \lambda^{1,0})$ can be
identified with the subspace of invariant classes in the sheaf cohomology space $H^1(M, {\cal O}
(T))$. By \cite[Proposition 4.1]{Sa} if $J$ is a smooth point of ${\cal C} ({\cg})$, then the tangent
space $T_J {\cal C}(\cg)$ to ${\cal C}(\cg)$ is contained in the complex subspace of $T_J {\cal C}$
(where ${\cal C} \cong  {GL (2n, \R)}/{GL (n, \C)}$ is the set of all almost complex 
structures on $\cg$) determined by $K = \ker \debar$.

By the above section we know that if $J_0$ is a rational complex structure, 
then $\HD {p} {q} (M) \cong
\HD {p}{q} (\cg^{\C})$. Given $J_0 \in {\cal C} ({\cg})$,  we know by 
\cite{KS} there exists a complete complex analytic
family $\{ M_t= (M, J_t) \mid J_t \in B \}$.

Let $\debar_t$ and $\Delta_t$ be respectively the  $\debar$-operator and the Laplacian 
determined by the global inner product   $g_t$ induced by an invariant Hermitian metric on $M$
 compatible with $J_t$. More precisely,
$$ 
\Delta_t = \debar_t^* \debar_t + \debar_t \debar_t^*,
$$
where $\debar_t^*$ is the adjoint of $\debar_t$ with respect to $g_t$. 

\medskip

\noindent{\bf Lemma 7}.  {\sl i) $\Delta_t$ sends $G$-invariant forms of 
type $(p, q)$ to $G$-invariants
forms. Moreover, the orthogonal complement with respect to $g_t$ of the invariant forms on the space
$\Lambda^{p,q}_t$ of $\Gamma$-invariant forms of type $(p,q)$ on $(M, J_t)$ is preserved by $\Delta_t$.

ii) $\HD{p}{q} (\cg^{\C})$ is a subspace of $\HD {p}{q} (M)$, for any 
invariant complex structure $J$ on $M$.}

\medskip

\noindent{\it Proof.} i) follows by the fact that $\debar_t$ and 
$\debar_t^*$ preserve $G$-invariant forms.

\noindent ii) We have to prove that there exists an injective homomorphism $$\pi: \HD {p}{q}(\cg^{\C})
 \to \HD {p}{q} (M).$$ By the decomposition
 $$
 \Lambda^{p,q} = {\cal H}^{p,q} \oplus Im \debar \oplus Im \debar^*,
 $$
 where ${\cal H}^{p,q}$ denotes the space of $\Gamma$-invariant harmonic forms of type $(p,q)$ on $M$,
 we have similarly, for the $G$-invariant forms, the decomposition
 $$
 \Lambda^{p,q}_{inv} = {\cal H}^{p,q}_{inv} \oplus Im \debar \vert_{inv} \oplus Im
\debar^*\vert_{inv},
$$
where $\HD {p}{q}(\cg^{\C}) \cong {\cal H}^{p,q}_{inv}$.

We can define $\pi (\omega)$ as the orthogonal projection of $\omega$ on 
${\cal H}^{p,q} = (Im \debar \oplus Im \debar^*)^{\perp}$, for any $\lbrack \omega \rbrack
\in \HD {p}{q}(\cg^{\C})$. To prove that $\pi$ is injective, suppose that $\pi (\omega) = 0$ on 
$\HD {p}{q} (M)$. Then $\pi (\omega) = \debar \varphi$. We may assume that 
$\varphi$ belongs to the
orthogonal complement to the $G$-invariant forms.
Since $\debar \varphi$ is $G$-invariant we have that also $\debar^* (\debar \varphi)$ is invariant
and then the inner product of $\varphi$ with $\debar^* (\debar \varphi)$ is zero and thus $\varphi$ must
be $G$-invariant, i.e. $\lbrack \omega \rbrack = 0$ in $\HD {p}{q} (\cg^{\C})$.

\medskip

\noindent{\bf Lemma 8}.  {\sl For any small deformation of the rational 
complex structure $J$ the isomorphism \rf{iso} holds.}

\medskip

\noindent{\it Proof.}
Using the same proof as in \cite[p.67-68]{KS} it is possible to show that for any $(p, q)$ and on any
$(M, J_t)$ there exists a complete orthonormal set of forms of type $(p,q)$ orthogonal to the
$G$-invariant ones $$
\{ e^{p,q}_{t1}, \ldots; f^{p,q}_{t1}, \ldots; g^{p,q}_{t1}, \ldots \} $$
such that 

(i) $\Delta_{t} e^{p,q}_{tj} = 0$, for any $j = 1, \ldots, d_q$, $d_q = \dim {\cal H}^{p,q}_{\perp}$
(where ${\cal H}^{p,q}_{\perp}$ is the space of harmonic forms of type $(p,q)$ orthogonal to the
$G$-invariant ones).

(ii) $\Delta_{t} f^{p,q}_{tj} = a_j^{p,q} (t) f^{p,q}_{tj}$;

(iii) $\Delta_{t} g^{p,q}_{tj} = a_j^{p,q -1} (t) g^{p,q}_{tj}$; with $a_j^{p,q} (t) > 0$.

Let $c$ be a positive constant and denote by $\nu^{p,q}(t)$ the number 
of eigenvalues $a_j^{p,q} (t)$ of $\Delta_t$ such that 
$a_j^{p,q} (t) < c$. Let $\Lambda^{p,q}_t (c)_{\perp}$ be the subspace spanned by the eigenforms of
$\Delta_t$ (orthogonal to the $G$-invariant ones) such that the corresponding eigenvalues are less
than $c$. Thus we have
$$
\dim \Lambda^{p,q}_t (c)_{\perp} = h_{p,q} (t)^{\perp} + \nu^{p,q}(t) + \nu^{p,q -1} (t),
$$
where $h_{p,q} (t)^{\perp}$ is the dimension of the orthogonal complement  $\HD {p}{q} (M_t)^{\perp}$
of the space of $G$-invariant forms $\HD {p}{q} (\cg^{\C})$ in $\HD {p}{q} (M_t)$.
Given a  point $J_{t_0} \in {\cal C} (\cg)$, we can choose $c$ such that 
$0 < c < a_j^{p,q} (t_0)$ for any $p, q = 0, \ldots, n$ and for any $j$.
 Moreover let $U$ be a sufficiently small neighbourhood of $t_0$ in $\cal C (\cg)$. Since any
eigenvalue $a_j^{p,q} (t)$ is a continuous function of $t$ \cite[Theorem 2, p.47]{KS} we have  that
the dimension of $\Lambda^{p,q}_t (c)_{\perp}$ is independent of $t \in U$ and  $\Lambda^{p,q}_{t_0}
(c)_{\perp} = \HD {p}{q} (M_{t_0})^{\perp}$. Thus we have
$$
h_{p,q} (t)^{\perp} + \nu^{p,q}(t) + \nu^{p,q -1} (t) = h_{p,q} (t_0)^{\perp}, 
$$
for any $t \in U$. This means that the function $h_{p,q} (t)^{\perp}$ 
is upper semicontinuous. Since  for any rational  complex structure  $J_0$, we have
that  $h_{p,q} (t_0)^{\perp} = 0$, so $h_{p,q} (t)^{\perp} = 0$ in some 
neighbourhood of $t_0$.
Consequently the set 
$$
\{ t \in {\cal C} (\cg) \mid h_{p,q} (t)^{\perp} = 0 \} = \{ t \in {\cal C} (\cg) \mid 
\HD{*}{*}(M_t) \cong \HD{*}{*}(\cg^{\C}) \}
$$
is a open set in any connected component of ${\cal C}(\cg)$ (with 
respect to the induced topology of ${\frak {gl}} (2n, \R)$ on ${\cal C}(\cg)$).

\section{Examples of compact nimanifolds with non rational
complex structures}\lb{Section finale}

Let $M = \Gamma\setminus G$ be the Iwasawa manifold. Recall that it can be constructed by taking as nilpotent Lie group $G$ the
complex Heisenberg group
$$
G = \left \{ \left( \begin{array} {ccc} 1&z_1&z_3\\ 0&1&z_2\\ 0&0&1 \end{array} \right) : z_i \in \C, i = 1, 2, 3 \right \},
$$
and, as lattice $\Gamma$, the subgroup of $G$ consisting of those matrices whose entries are Gaussian integers. It is known
that the 1-forms
${\omega}_1 = d z^1, {\omega}_2 = d z_2, {\omega}_3  = d z_3 + z_1 d z_2,$ are left invariant on $G$. $G$ has structure
equations
$$
d \omega_1 = d \omega_2 = 0, \qquad d \omega_3 = \omega_1 \wedge \omega_2.
$$
If one regards $G$ as a real Lie group and sets 
$$
\omega_1 =: e^1 + i e^2, \quad \omega_2 =: e^3 + i e^4, \quad \omega_3 =: e^5 + i e^6 \eqno(2)
$$
then $( e^i)$ is a real basis of ${\frak g}^*$ such that:
$$
\left \{ \begin{array} {l} d e^i = 0, \qquad \qquad \qquad \qquad 1 \leq i \leq 4,\\
de^5= e^1 \wedge e^3 + e^4 \wedge e^2,\\
d e^6 = e^1 \wedge e^4 + e^2 \wedge e^3. 
\end{array} \right.
$$
By \cite{Na},  the 
bi-invariant complex structure $J_0$ defined by (2) on the complex Heisenberg group has a deformation space of positive
dimension. Then the complex structure associated to the subspace 
$$
<\omega_1,  \omega_2, \omega_3 + t {\bar \omega}_1>
$$
of ${\frak g}^*_{\C}$ belongs to the  orbit of $J_0$ in ${\cal C} ({\frak g})$ with respect the group of automorphisms of
$\frak g$ (see
\cite[Section 4]{Sa}) and it  is  not rational for appropriate
$t\in \C$. 

In the same way one can see that every compact nilmanifold $M=G/\Gamma$ of real dimension $2n$, with at least one (rational)
complex structure
$J_0$ having a non trivial deformation, has a non rational complex structure. Indeed, by
\cite[Theorem 1.3]{Sa}, given the complex structure  $J_0$ it
is possible to construct a basis of left invariant $(1,0)$-forms $\{ \omega_1, \ldots, \omega_n \}$ such that $d
\omega_{i + 1}$ belongs to the ideal generated by the set $\{ \omega_1, \ldots, \omega_i \}$ ($i=0, \ldots, n-1$) in the
complexified exterior algebra.  Then the complex structure associated to the subspace 
$$
<\omega_1, \ldots, \omega_n + t \bar \omega_1>
$$
is a not rational complex structure on $M$ for  appropriate
$t\in \C$.

\bigskip \bigskip

\renewcommand{\thebibliography}{\list{[\arabic{enumi}]\hfil}
{\settowidth\labelwidth{18pt}\leftmargin\labelwidth\advance
\leftmargin\labelsep\usecounter{enumi}}\def\newblock{\hskip.05em}
\sloppy\sfcode`\.=1000\relax}\newcommand{\bi}{\vspace{-3pt}\bibitem}

\end{document}